\documentclass[12pt]{article}
\usepackage{setspace}
\usepackage[margin=1in]{geometry}
\usepackage{amsmath}
\usepackage{amssymb}
\usepackage{graphicx,multirow}
\usepackage{algorithm}
\usepackage{algorithmic}
\usepackage{subfigure}
\usepackage{multicol,float}
\usepackage{bbm}
\usepackage{url}
\usepackage{tikz,pgfplots}
\usepackage{amsthm}
\usepackage{authblk}
\usepackage[round,sort]{natbib}
\usepackage[pdfborder={0 0 0}]{hyperref}
\def\EE{\mathbb{E}}
\def\PP{\mathbb{P}}

\begin{document}
\doublespacing

\title{Kelly betting on horse races with uncertainty\\in probability estimates}

\author{Michael R. Metel}

\affil{Laboratoire de Recherche en Informatique, Universit\'e Paris-Sud, Orsay, France\\
\url{metel@lri.fr}}

\maketitle

\begin{abstract}
\doublespacing
We investigate the problem of gambling with uncertainty in outcome probabilities. Stochastic optimization models are proposed for optimal investing on events with mutually exclusive outcomes when probabilities are estimated using multinomial logistic regression. Special attention is given to the case of there being two outcomes, and the general case of many outcomes. An empirical study using simulated data was conducted where the loss of return from probability estimation error is observed, and superior returns are achieved taking it into consideration.
\end{abstract}

\section{Introduction}

The \cite{kelly1956} criterion is a powerful tool for decision making in the world of gambling and investing, answering the question of how much to wager by maximizing the asymptotic exponential rate of return. One limitation is that exact knowledge of outcome probabilities and payouts is assumed, which in general is not available for opportunities which have the potential for profitability, such as in sports betting and the stock market. Replacing the actual values with estimates leads to overbetting \citep{maclean1992}, resulting in higher risk with lower returns. We might assume unbiased errors will cancel out over time, but when using the Kelly criterion, outcomes with overestimated values will consistently look more favourable than in reality, with larger wagers placed on them than should be. The most popular way to mitigate this problem is through the use of a fractional Kelly strategy, which wagers a fixed fraction of the amount prescribed by the Kelly criterion, which has been shown to possess favourable risk-return properties by \cite{maclean1992}. This can be seen as a conservative strategy, where given an estimated upper bound of possible over betting, we correct by dividing all wagers by this amount.
This technique has been endorsed and successfully used in practice by people such as \cite{ben08} and \cite{thorp2006}, with betting half the Kelly amount being popular amongst gamblers \citep{Pound05}.
More recently, research has been done which directly considers the inherent uncertainty in probability estimates. \cite{Baker13} derived a shrinkage factor for two outcome gambling settings, which gives the optimal fraction of the Kelly amount to wager.\\

In this paper we develop methodologies for Kelly style betting on events with many possible outcomes taking into account the uncertainty in probability estimates. For clarity we present our work using the example of horse race betting, though the ideas can be applied more broadly in domains beyond sports betting, and should be of interest to people concerned with decision making under uncertainty in general. A significant amount of academic research has been done on horse racing, see \citep{hau81}, largely due to the fact that it can be seen as a microcosm of financial exchanges, with findings in this arena having wider implications in finance, economics, and decision theory.

\section{Optimal investment problem}

We are interested in the following problem, which we will describe in the setting of optimal wagering on a horse race.

\begin{alignat}{6}
&\max&&\text{ }\sum_{h=1}^n\pi_{h}\log(x_{h}O_{h}+w-\sum_{i=1}^nx_i)\label{eq:SF}\tag{P}\\
&\mbox{s.t. }&&\sum_{h=1}^nx_h\leq w \nonumber\\
&&&x_{h}\geq 0\hspace{50 pt}h=1,...,n \nonumber
\end{alignat}

The objective is to maximize the Kelly criterion, which is equivalent to maximizing the expected logarithm of wealth. There are $n$ horses in the race, and the probability of horse $h$ winning is $\pi_h$. $x_h$ is the amount we will wager on horse $h$, $O_h$ is the decimal payout odds for horse $h$, which we assume are fixed, and $w$ is our current wealth. The constraints ensure our wagers are sound. Exact optimal solutions can be found for this problem using the algorithm devised by \cite{smoc2010}.\\

\section{Outcome probability estimation}

\label{sec:OP}

The most popular method of estimating horse racing win probabilities is by multinomial logistic regression, which was first proposed by \cite{Bolt86}.
We create a linear predictor function to value each horse $h$, $\beta'v_{h}$, where $v_{h}$ is the vector of factor values and $\beta'$ is the transposed vector of regression coefficients. Each horse is then assigned winning probabilities $\pi_{h}=\frac{e^{\beta'v_{h}}}{\sum_{i=1}^{n}e^{\beta'v_{i}}}$. Given $R$ historical races, the log-likelihood function is
$$\ln L(\beta)=\sum_{r=1}^R \beta'v_{w^r}-\ln \sum_{i=1}^{n_r}e^{\beta'v_{i^r}},$$
where $w^r$ denotes the index of the winning horse in race $r$ and $n_r$ is the number of horses. This function is
concave~\cite[p. 72]{boyd2004}, so we can find a maximum likelihood estimate $\hat{\beta}$ using standard unconstrained optimization algorithms. The $j^{th}$ component of the score function, or the gradient of $\ln L(\beta)$ is
$$\frac{\partial\ln L(\beta)}{\partial \beta_j}=\sum_{r=1}^R v_{w^r_j}-\frac{\sum_{i=1}^{n_r}v_{i^r_j}e^{\beta'v_{i^r}}}{\sum_{i=1}^{n_r}e^{\beta'v_{i^r}}}.$$
The $jk$ cell of the curvature, or the negative Hessian is
$$\frac{-\partial^2\ln L(\beta)}{\partial\beta_j\partial\beta_k}=\sum_{r=1}^R \frac{(\sum_{i=1}^{n_r}v_{i^r_j}v_{i^r_k}e^{\beta'v_{i^r}})(\sum_{i=1}^{n_r}e^{\beta'v_{i^r}})-
(\sum_{i=1}^{n_r}v_{i^r_k}e^{\beta'v_{i^r}})(\sum_{i=1}^{n_r}v_{i^r_j}e^{\beta'v_{i^r}})}{(\sum_{i=1}^{n_r}e^{\beta'v_{i^r}})^2}.$$
The observed Fisher information, $I(\hat{\beta})$, is the curvature at $\hat{\beta}$. The maximum likelihood estimate $\hat{\beta}$ is consistent and asymptotically normal with covariance $I^{-1}(\hat{\beta})$ \citep{mcfad74}.\\

There are too many potential factors contributing to the outcome of a horse race to believe we have the actual model of outcome probabilities, so we must consider model misspecification. The variance of the score terms is
$V(\hat{\beta})=\sum_{r=1}^R \frac{\partial\ln L_r(\hat{\beta})}{\partial \beta}\frac{\partial\ln L_r(\hat{\beta})}{\partial \beta}'$, where $\frac{\partial\ln L_r(\beta)}{\partial \beta}$ is the gradient of the log-likelihood function of the $r^{th}$ race. We can then calculate the sandwich estimate of the covariance matrix of $\hat{\beta}$ as $\Sigma=I(\hat{\beta})^{-1}V(\hat{\beta})I(\hat{\beta})^{-1}$, which is robust to our misspecification, and model our parameters as $\beta\sim N(\hat{\beta},\Sigma)$. Let
$\pi^s_{h}=\frac{e^{\EE(v_{h})}}{\sum_{i=1}^{n}e^{\EE(v_i)}}=\frac{e^{\hat{\beta}'v_{h}}}{\sum_{i=1}^{n}e^{\hat{\beta}'v_{i}}}$,
which is the maximum likelihood probability estimate typically used after completing a regression analysis. We denote (\ref{eq:SF}) with $\pi_h=\pi^s_h$ as (S) and call this the standard model, which will be used to compare other models to.

\section{Optimization models considering uncertainty in probability estimates}

We now attempt to move beyond the standard model (S) by considering the uncertainty of our probability estimates using different techniques from stochastic optimization. A natural approach would be to maximize the expected value of our objective,
$$\EE\sum_{h=1}^n \pi_h(\beta)\log(x_{h}O_h+w-\sum_{i=1}^nx_i),$$
resulting in (\ref{eq:SF}) with probabilites $\EE(\pi_h)=\EE\frac{e^{\beta'v_{h}}}{\sum_{i=1}^{n}e^{\beta'v_{i}}}$.
As random variables our win probabilities,
$\pi_{h}=\frac{e^{\beta'v_{h}}}{\sum_{i=1}^{n}e^{\beta'v_{i}}}$, follow a logistic-normal distribution \citep{atch1980} for which $\EE(\pi_h)$, to the best of our knowledge, is not representable in a simple form. $\EE(\pi_h)$ can be estimated using Monte Carlo integration by generating $N$ random samples of $\beta$ and taking our expected win probabilites as $\pi^{mc}_h=\frac{1}{N}\sum_{i=1}^N\frac{e^{(\beta^i)'v_{h}}}{\sum_{j=1}^{n}e^{(\beta^i)'v_{j}}}$. We denote this formulation as (Emc).\\

We can also find a lower bound using the normal moment generating function,
$$\EE\left(\frac{1}{\pi_h}\right)=\EE\left(\sum_{i=1}^{n}e^{\beta'(v_i-v_h)}\right)=\frac{\sum_{i=1}^{n}e^{\hat{\beta}'v_i+\frac{1}{2}(v_i-v_h)'\Sigma(v_i-v_h)}}{e^{\hat{\beta}'v_h}},$$ and from Jensen's inequality, $\EE(\pi_h)\geq \frac{1}{\EE(\frac{1}{\pi_h})}$, giving us the lower bound
$$\EE^{lb}(\pi_h)=\frac{e^{\hat{\beta}'v_h}}{\sum_{i=1}^{n}e^{\hat{\beta}'v_i+\frac{1}{2}(v_i-v_h)'\Sigma(v_i-v_h)}}$$
for each outcome probability. Using these lower bounds gives us a conservative estimate of outcomes. The sum of  probabilities $\EE^{lb}(\pi_h)$ will in general not equal 1, so when betting on an outcome $h$, this formulation underweights the event of $h$ not occurring. To overcome this we have the following formulation which contains an extra outcome with probability $1-\sum_{i=1}^n \EE^{lb}(\pi_i)$ where all money wagered is lost.

\begin{alignat}{6}
&\max&&\text{ }\sum_{h=1}^n\EE^{lb}(\pi_h)\log(x_{h}O_{h}+w-\sum_{i=1}^nx_i)+\left(1-\sum_{h=1}^n\EE^{lb}(\pi_h)\right)\log(w-\sum_{i=1}^nx_i) \label{eq:Elb}\tag{Elb}\\
&\mbox{s.t. }&&\sum_{h=1}^nx_h\leq w \nonumber\\
&&&x_{h}\geq 0\hspace{50 pt}h=1,...,n \nonumber
\end{alignat}

As we want to limit overbetting from probability estimation error, we can add solution robustness to probability estimation. The optimization problem is rewritten so that the uncertainty is in the constraints and then we ensure a minimal objective value holds for a given probability through the use of a chance constraint.
\begin{alignat}{6}
&\max&&\text{ }t\label{eq:RO}\tag{CC}\\
&\mbox{s.t. }&&\PP(t\leq\sum_{h=1}^n\pi_h\log(x_{h}O_h+w-\sum_{i=1}^nx_i))\geq 1-\alpha\nonumber\\
&&&\sum_{h=1}^nx_h\leq w \nonumber\\
&&&x_{h}\geq 0\hspace{50 pt}h=1,...,n \nonumber
\end{alignat}

where $\alpha<0.5$. When using (S), we are only optimizing over a point estimate of $\pi$. Assuming this estimate differs from the actual value of $\pi$, we could very well be placing a wager with actual expected log wealth lower than our current $\log(w)$. With this chance constraint we can ensure that our solution generates a positive expected return over a high proportion of possible values of $\pi$. Further, we can choose a solution $x$ which will generate the highest return over $(1-\alpha)\%$ of potential values of $\pi$, avoiding large misplaced wagers.\\

In the case where there are only two outcomes, we can achieve the exact solution by solving the following optimization program.
\begin{alignat}{6}
&\max&&\text{ }t\label{eq:CC2E}\tag{CC2}\\
&\mbox{s.t. }&&t\leq \pi^H_1\log(x_{1}(O_1-1)+w-x_2)+\pi^L_2\log(x_{2}(O_2-1)+w-x_1)\nonumber\\
&&&t\leq \pi^L_1\log(x_{1}(O_1-1)+w-x_2)+\pi^H_2\log(x_{2}(O_2-1)+w-x_1)\nonumber\\
&&&x_1+x_2\leq w \nonumber\\
&&&x_1,x_2\geq 0\nonumber
\end{alignat}
where
\begin{alignat}{6}
&\pi^H_1=\frac{e^{\Phi^{-1}(1-\alpha)\sigma+\hat{\beta}'v_1}}{e^{\Phi^{-1}(1-\alpha)\sigma+\hat{\beta}'v_1}+e^{\hat{\beta}'v_2}},&\pi^L_2=\frac{e^{\hat{\beta}'v_2}}{e^{\Phi^{-1}(1-\alpha)\sigma+\hat{\beta}'v_1}+e^{\hat{\beta}'v_2}}\nonumber\\
&\pi^L_1=\frac{e^{\hat{\beta}'v_1}}{e^{\hat{\beta}'v_1}+e^{\Phi^{-1}(1-\alpha)\sigma+\hat{\beta}'v_2}},&\pi^H_2=\frac{e^{\Phi^{-1}(1-\alpha)\sigma+\hat{\beta}'v_2}}{e^{\hat{\beta}'v_1}+e^{\Phi^{-1}(1-\alpha)\sigma+\hat{\beta}'v_2}}\nonumber
\end{alignat}

and $\sigma^2=(v_1-v_2)'\Sigma(v_1-v_2)$. Intuitively, our probabilities $\pi_1$ and $\pi_2$ are reweighted depending on which outcome is more favourable, where the first constraint puts more weight on outcome one, which will be tight when outcome two is more favourable, and the second constraint puts more weight on outcome two, for when one is more favourable. A derivation of this program can be found in the appendix in the subsection {\it Derivation of (CC2)}.\\

We now focus on the case of more than two outcomes. The following optimization problem is an approximation of \ref{eq:RO}, where we have taken $S$ iid samples $\pi^s$ of the outcome probabilities and want to satisfy the chance constraint over this empirical distribution by ensuring the constraint $t\leq\sum_{h=1}^n\pi^s_h\log(x_{h}O_h+w-\sum_{i=1}^nx_i)$ is not satisfied over no more than  $S\alpha$ samples. $M$ is chosen sufficiently large so as not to restrict the value of $t$. The convergence of the optimal objective value and solution set of this approximation to (\ref{eq:RO}) in the limit is established in \cite[p. 211]{shap2009}.

\begin{alignat}{6}
&\max&&\text{ }t\label{eq:CCN}\tag{CCN}\\
&\mbox{s.t. }&&t\leq \sum_{h=1}^n\pi^s_h\log(x_{h}O_h+w-\sum_{i=1}^nx_i)+z_sM\hspace{50 pt}s=1,...,S\nonumber\\
&&&\sum_{s=1}^S z_s\leq S\alpha\nonumber\\
&&&\sum_{h=1}^nx_h\leq w \nonumber\\
&&&x_{h}\geq 0\hspace{70 pt}h=1,...,n \nonumber\\
&&&z_s\in\{0,1\}\hspace{50 pt}s=1,...,S \nonumber
\end{alignat}

This problem is quite challenging, and is only practical for very small sample choices of $S$. We use a simple heuristic to find the $S\alpha$ worst constraints, setting their binary values to $1$, and then proceed to solve the now convex problem with decision variables $x_h$ and $t$. We first solve (\ref{eq:CCN}) with  $\alpha=0$, or for $z_s=0$ $\forall s\in S$, sort $t-\sum_{h=1}^n\pi^s_h\log(x_{h}O_h+w-\sum_{i=1}^nx_i)$ in descending order, set $z_s=1$ for the first $S\alpha$ corresponding constraints, and then resolve (\ref{eq:CCN}) with fixed $z_s$.\\

We consider a final model, combining the previous two, where we are maximizing the expectation of log wealth subject to a chance constraint which ensures our solution does not have a negative true expected return with high probability.

\begin{alignat}{6}
&\max&&\text{ }\EE \sum_{h=1}^n\pi_h\log(x_{h}O_h+w-\sum_{i=1}^nx_i)\label{eq:ECC}\tag{ECC}\\
&\mbox{s.t. }&&\PP(\log(w)\leq\sum_{h=1}^n\pi_h\log(x_{h}O_h+w-\sum_{i=1}^nx_i))\geq 1-\alpha\nonumber\\
&&&\sum_{h=1}^nx_h\leq w \nonumber\\
&&&x_{h}\geq 0\hspace{50 pt}h=1,...,n \nonumber
\end{alignat}

The implementation uses the objective of (\ref{eq:SF}) with probabilities $\pi^{mc}_h$ and the constraint set of either (\ref{eq:CC2E}) or (\ref{eq:CCN}), with $t$ replaced with $\log(w)$. We label these formulations (ECC2) and (ECCN) respectively.

\section{Empirical model comparison}

\label{sec:ES}

We seek to compare the performance of standard Kelly betting, fractional Kelly betting and Kelly betting considering the uncertainty in probability estimates, using simulated data to allow for accurate testing of the different methodologies. We simulated the distribution of $\beta$ by taking $\hat{\beta}$ as a standard normal random vector of size $m=10$. For each component $\hat{\beta}_i$, a corresponding standard deviation $\sigma_i$ was simulated between $[0,\sigma'_i]$, where $\sigma'_i=\frac{-|\beta_i|}{\Phi^{-1}(0.025)}$. This implies that the p-value of $\beta_i$ is not greater than $0.05$ if the data were the result of a regression analysis. $\Sigma$ was then taken as a diagonal matrix consisting of $\sigma^2$. The true values of $\beta$, $\beta^t$, corresponding to the true outcome probabilities were simulated by taking a random sample from the distribution of $\beta$. We generated an $m \times n$ matrix $F$ of standard uniform random variables representing the $m$ factor values of the $n$ horses. We then calculated the true outcome probabilities as

$$\pi^t_h=\frac{e^{{\beta^t}' F_{h}}}{\sum_{i=1}^ne^{{\beta^t}'F_{i}}}$$

In order to limit the variation of the empirical testing, we compared models not based on their final simulated wealth, but by their expected exponential return over the true probability distribution $\pi^t$. Given our optimal solution $x^*$ using any technique in race $i$, the expected exponential return is calculated as

$$\EE\log\left(\frac{w_{i}}{w_{i-1}}\right)=\sum_{h=1}^n\pi^t_h\log\left(\frac{x^*_{h}O_h+w_{i-1}-\sum_{i=1}^nx^*_i}{w_{i-1}}\right)$$

We can then calculate the expected return over the entire sample of size $T$ as

$$\EE\log\left(\frac{w_{T}}{w_{0}}\right)=\sum_{i=1}^T\EE\log\left(\frac{w_{i}}{w_{i-1}}\right)$$

This should reduce the variance in the result, as our performance statistic does not depend on the outcomes of the races. To further limit variance, we use a fixed, identical payout odd for each outcome.\\

Four experiments were conducted, each consisting of 2500 trials, with details in Table \ref{T1}.

\begin{table}[H]
\centerline{
\resizebox{0.3\textwidth}{!}{
\renewcommand{\arraystretch}{1}
\begin{tabular}{lrr}
\bf{Experiment}&\bf{n}&$\mathbf{O_h}$\\
  \hline
(E1)    &2&1.1\\
(E2)    &2&1.2\\
(E3)    &10&2\\
(E4)    &30&4\\
\hline
\end{tabular}}}
\caption{Experiment details} \label{T1}
\end{table}

All chance constraints were tested using $\alpha=0.4$, $0.25$, and $0.1$.
$\pi^{mc}_h$ was estimated using $1,000,000$ samples when $n=2$, and $2,000,000$ when $n=10$ and $n=30$.
The chance constraints were estimated using $1,000$ samples for $n=10$ and $2,000$ samples for $n=30$.
All experiments were done on a Windows 10 Home 64-bit, Intel Core i5-7200U 2.5GHz processor with 8 GB of RAM, in Matlab R2017a using {\it fmincon}.\\

Result data is presented in Table \ref{T5} for all experiments.
(T) is the best that could be achieved in the experiment, using the true probabilities $\pi^t_h$ in (\ref{eq:SF}). (F) is using a 50\% fractional Kelly strategy of (S). (CCx) stands for (CC2) or (CCN) depending on the experiment, and likewise for (ECCx).
The numbers under each experiment name are the total expected return for each model, with the final column being the sum over all experiments. We see in total, (\ref{eq:Elb}), (Emc), and (CCx) and (ECCx) with $\alpha=0.4$ outperformed (\ref{eq:SF}), with (Emc) and (ECCx) with $\alpha=0.4$ outperforming (S) in all experiments, and (ECCx) with $\alpha=0.4$ performing best overall. Though the best result was achieved using a chance constraint, we can see that in general its use as a risk measure against uncertainty in probability estimation is overly aggressive for smaller values of $\alpha$, dampening long term growth.\\

\begin{table}[H]
\centerline{
\resizebox{0.65\textwidth}{!}{
\begin{tabular}{llrrrrr}
   \bf{Model}&&\bf{(E1)}&\bf{(E2)}&\bf{(E3)}&\bf{(E4)}&\bf{Sum}\\
  \hline
(T)    &             &3.051&16.861&4.385&3.166&27.463\\
(S)    &             &1.583&12.214&2.383&1.954&18.134\\
(F)    &             &1.210&8.754&1.788&1.515&13.268\\
(Elb)  &             &1.852&12.179&2.353&1.823&18.206\\
(Emc)  &             &1.820&12.297&2.390&1.964&18.471\\
(CCx) &$ \alpha=0.40$  &1.786&12.288&2.345&1.834&18.253\\
(CCx) &$ \alpha=0.25$ &1.694&11.371&1.964&1.426&16.454\\
(CCx) &$ \alpha=0.10$ &1.080&9.027&1.206&0.846&12.158\\
(ECCx)&$ \alpha=0.40$  &1.825&12.288&2.390&1.980&18.484\\
(ECCx)&$ \alpha=0.25$ &1.834&11.983&2.319&1.728&17.863\\
(ECCx)&$ \alpha=0.10$ &1.433&10.355&1.647&1.164&14.599\\
\hline
\end{tabular}}}
\caption{Expected total returns} \label{T5}
\end{table}

\section{Conclusion and future research}
\label{s:ME}

We have investigated different stochastic optimization models for Kelly style betting on mutually exclusive outcomes considering probability estimation uncertainty stemming from multinomial logistic regression. An empirical study using simulated data was conducted to compare performance. The large difference in long term growth when using the true probability outcomes versus relying on estimates in our experiments display the significance of probability estimation error in decision making, and the challenges for those attempting to maximize return in speculative markets. Improvements in long term growth have been found, first by considering the uncertainty in outcome probabilities when calculating the expected log wealth, and with a mild use of a chance constraint, which will likely need to be calibrated in each application to find the proper balance of preventing losses from uncertainty without overly dampening the potential to capture positive returns.\\

The presentation of the material in this paper has focused on the application of betting on horse racing, but the ideas are applicable to a general investment setting. New research adapting the methods presented to more general return settings beyond multinomial logistic regression uncertainty and mutual exclusive events would be interesting, with applications such as investing in a portfolio of stocks following geometric Brownian motions considering parameter uncertainty.

\section*{Acknowledgements}

The author thanks the anonymous referees for their valuable comments. This work was supported by the Digiteo Chair C\&O program.

\bibliographystyle{plainnat}
\bibliography{DARef_FINAL}

\begin{thebibliography}{13}
\providecommand{\natexlab}[1]{#1}
\providecommand{\url}[1]{\texttt{#1}}
\expandafter\ifx\csname urlstyle\endcsname\relax
  \providecommand{\doi}[1]{doi: #1}\else
  \providecommand{\doi}{doi: \begingroup \urlstyle{rm}\Url}\fi

\bibitem[Atchison and Shen(1980)]{atch1980}
J~Atchison and SM~Shen.
\newblock Logistic-normal distributions: {S}ome properties and uses.
\newblock \emph{Biometrika}, 67\penalty0 (2):\penalty0 261--272, 1980.

\bibitem[Baker and McHale(2013)]{Baker13}
RD~Baker and IG~McHale.
\newblock Optimal {B}etting {U}nder {P}arameter {U}ncertainty: Improving the
  {K}elly {C}riterion.
\newblock \emph{Decision Analysis}, 10\penalty0 (3):\penalty0 189--199, 2013.

\bibitem[Benter(1994)]{ben08}
W~Benter.
\newblock {C}omputer {B}ased {H}orse {R}ace {H}andicapping and {W}agering
  {S}ystems: A {R}eport.
\newblock In DB~Hausch, VSY Lo, and WT~Ziemba, editors, \emph{Efficiency of
  {R}acetrack {B}etting {M}arkets}, pages 183--198. World Scientific, 1994.

\bibitem[Bolton and Chapman(1986)]{Bolt86}
RN~Bolton and RG~Chapman.
\newblock Searching for {P}ositive {R}eturns at the {T}rack: A {M}ultinomial
  {L}ogic {M}odel for {H}andicapping {H}orse {R}aces.
\newblock \emph{Management Science}, 32\penalty0 (8):\penalty0 1040--1060,
  1986.
\newblock ISSN 0025-1909.

\bibitem[Boyd and Vandenberghe(2004)]{boyd2004}
S~Boyd and L~Vandenberghe.
\newblock \emph{Convex {O}ptimization}.
\newblock Cambridge {U}niversity {P}ress, 2004.

\bibitem[Hausch et~al.(1981)Hausch, Ziemba, and Rubinstein]{hau81}
DB~Hausch, WT~Ziemba, and M~Rubinstein.
\newblock Efficiency of the {M}arket for {R}acetrack {B}etting.
\newblock \emph{Management Science}, 27\penalty0 (12):\penalty0 1435--1452,
  1981.

\bibitem[Kelly(1956)]{kelly1956}
JL~Kelly.
\newblock A {N}ew {I}nterpretation of {I}nformation {R}ate.
\newblock \emph{Information Theory, IRE Transactions on}, 2\penalty0
  (3):\penalty0 185--189, 1956.

\bibitem[MacLean et~al.(1992)MacLean, Ziemba, and Blazenko]{maclean1992}
LC~MacLean, WT~Ziemba, and G~Blazenko.
\newblock Growth {V}ersus {S}ecurity in {D}ynamic {I}nvestment {A}nalysis.
\newblock \emph{Management Science}, 38\penalty0 (11):\penalty0 1562--1585,
  1992.

\bibitem[McFadden(1974)]{mcfad74}
D~McFadden.
\newblock Conditional {L}ogit {A}nalysis of {Q}ualitative {C}hoice {B}ehavior.
\newblock In Zarembka P, editor, \emph{Frontiers in {E}conometrics}, pages
  105--142. Academic Press, New York, 1974.

\bibitem[Poundstone(2005)]{Pound05}
W~Poundstone.
\newblock \emph{Fortune's {F}ormula: {T}he {U}ntold {S}tory of the {S}cientific
  {B}etting {S}ystem that {B}eat the {C}asinos and {W}all {S}treet}.
\newblock Hill and Wang, 2005.

\bibitem[Shapiro et~al.(2009)Shapiro, Dentcheva, and Ruszczy{\'n}ski]{shap2009}
A~Shapiro, D~Dentcheva, and A~Ruszczy{\'n}ski.
\newblock \emph{Lectures on stochastic programming: modeling and theory}.
\newblock SIAM, 2009.

\bibitem[Smoczynski and Tomkins(2010)]{smoc2010}
P~Smoczynski and D~Tomkins.
\newblock An explicit solution to the problem of optimizing the allocations of
  a bettor's wealth when wagering on horse races.
\newblock \emph{Mathematical Scientist}, 35\penalty0 (1):\penalty0 10--17,
  2010.

\bibitem[Thorp(2006)]{thorp2006}
EO~Thorp.
\newblock The {K}elly {C}riterion in {B}lackjack, {S}ports {B}etting, and the
  {S}tock {M}arket.
\newblock In SA~Zenios and WT~Ziemba, editors, \emph{Handbook of Asset and
  Liability Management, Volume I}, pages 385--428. Elsevier, 2006.

\end{thebibliography}

\section*{Appendix}

\subsection*{Derivation of (CC2)}

We need to find an equivalent deterministic constraint for the chance constraint $\PP(t\leq \pi_1\log(x_{1}(O_1-1)+w-x_2)
+\pi_2\log(x_{2}(O_2-1)+w-x_1))\geq 1-\alpha$. For simplicity, let $W_1=x_{1}(O_1-1)+w-x_2$ and
$W_2=x_{2}(O_2-1)+w-x_1$, then the chance constraint can be written as
\begin{alignat}{6}
\PP(t\leq \pi_1\log(W_1)+(1-\pi_1)\log(W_2))\geq 1-\alpha.\label{eq:CC2}
\end{alignat}

Given that $\pi_1\in(0,1)$, the only time $t\geq\max(\log(W_1),\log(W_2))$ is feasible in (\ref{eq:CC2}) and (\ref{eq:CC2E}) is when $t=\log(W_1)=\log(W_2)$, with all other instances being infeasible in both. The case where $t\leq\min(\log(W_1),\log(W_2))$ is always feasible in (\ref{eq:CC2}) and (\ref{eq:CC2E}). We now consider the case where $\log(W_2)<t<\log(W_1)$. Rearranging and taking $\pi_1=\frac{e^{\beta'v_1}}{e^{\beta'v_1}+e^{\beta'v_2}}$, the chance constraint equals
\begin{alignat}{6}
&\PP(t-\log(W_2)\leq \pi_1(\log(W_1)-\log(W_2))\geq 1-\alpha\nonumber\\
&\PP(e^{\beta'(v_2-v_1)}\leq \frac{\log(W_1)-\log(W_2)}{t-\log(W_2)}-1)\geq 1-\alpha\nonumber\\
&\PP(\beta'(v_2-v_1)\leq \log\left(\frac{\log(W_1)-t}{t-\log(W_2)}\right))\geq 1-\alpha\nonumber\\
&\log\left(\frac{\log(W_1)-t}{t-\log(W_2)}\right)\geq \Phi^{-1}(1-\alpha)\sigma+\mu\nonumber
\end{alignat}

where $\mu=\hat{\beta}'(v_2-v_1)$ and  $\sigma^2=(v_1-v_2)'\Sigma(v_1-v_2)$. Rearranging,
\begin{alignat}{6}
&\frac{1}{1+e^{\Phi^{-1}(1-\alpha)\sigma+\mu}}\log(W_1)+\frac{e^{\Phi^{-1}(1-\alpha)\sigma+\mu_2}}{1+e^{\Phi^{-1}(1-\alpha)\sigma+\mu}}\log(W_2)\geq t\nonumber\\
&\frac{e^{\hat{\beta}'v_1}}{e^{\hat{\beta}'v_1}+e^{\Phi^{-1}(1-\alpha)\sigma+\hat{\beta}'v_2}}\log(W_1)+\frac{e^{\Phi^{-1}(1-\alpha)\sigma+\hat{\beta}'v_2}}{e^{\hat{\beta}'v_1}+e^{\Phi^{-1}(1-\alpha)\sigma+\hat{\beta}'v_2}}\log(W_2)\geq t\nonumber
\end{alignat}

which is the second constraint in (\ref{eq:CC2E}). In the case where $\log(W_1)<t<\log(W_2)$, we will get the same result, but with the probability shifted towards outcome 1, or the first constraint in (\ref{eq:CC2E}). Finally, we must show that the correct constraint will be active. When $\log(W_2)<t<\log(W_1)$,

\begin{alignat}{6}
&(\pi^H_1-\pi^L_1)\log(W_2)\leq (\pi^H_1-\pi^L_1)\log(W_1)\nonumber\\
&\pi^L_1(\log(W_1)-\log(W_2))\leq \pi^H_1(\log(W_1)-\log(W_2))\nonumber\\
&\pi^L_1\log(W_1)-(1-\pi^H_2)\log(W_2)\leq \pi^H_1\log(W_1)-(1-\pi^L_2)\log(W_2)\nonumber\\
&\pi^L_1\log(W_1)+\pi^H_2\log(W_2)\leq \pi^H_1\log(W_1)+\pi^L_2\log(W_2),\nonumber
\end{alignat}

and so the right hand side of the second constraint is not greater than the right hand side of the first constraint in (\ref{eq:CC2E}), implying the second constraint will be active. We can similarly show in the case when $\log(W_1)<t<\log(W_2)$ the opposite holds, that $\pi^H_1\log(W_1)+\pi^L_2\log(W_2)\leq  \pi^L_1\log(W_1)+\pi^H_2\log(W_2)$, implying the first constraint will be active as desired.

\end{document}